\documentclass[12pt]{article}
\usepackage{amsmath,amssymb,amsthm,array,amsfonts}
\usepackage{mathrsfs}   
\newcommand{\nmr}{\symbol{"9D}}

\newcommand{\R}[1]{\mathscr{R}(#1)}
\newcommand{\N}[1]{\mathscr{N}(#1)}
\newcommand{\ov}[1]{\mathscr D(#1)}

\newcommand{\e}{\varepsilon}

\newcommand{\oD}{\mathcal D}
\newcommand{\oDs}{\mathcal D^*}

\newcommand{\cl}{\mathrm cl\,}

\newcommand{\dt}{\mathrm dt}
\newcommand{\ds}{\mathrm ds}
\newcommand{\dF}{\frac d{dt}F}
\newcommand{\dFt}{\frac d{dt}F'}
\newcommand{\dfn}{\,:=\,}

\newcommand{\f}[1][\cdot]{f(#1)}

\newcommand{\ft}{\f[t]}
\newcommand{\z}[1][\cdot]{z(#1)}

\newcommand{\x}[1][\cdot]{x(#1)}
\newcommand{\xt}{x(t)}
\newcommand{\hx}[1][\cdot]{\hat{x}(#1)}

\renewcommand{\v}[1][\cdot]{v(#1)}

\newcommand{\Rm}{\mathbb{R}^m}
\newcommand{\Wn}{\mathbb{W}_2^n}

\newcommand{\WnF}{\mathbb{W}_{2}^\mathrm{F}}
\newcommand{\WmF}{\mathbb{W}_{2}^\mathrm{F'}}

\newcommand{\Ln}{\mathbb{L}_2^n}
\newcommand{\Lm}{\mathbb{L}_2^m}

\newcommand{\LM}{\mathbb{L}_2([t_0,T],\Rm)}

\newtheorem{thm}{Theorem}

\newtheorem{rem}{Remark}

\begin{document}
\begin{center}
\textbf{ON SOME PROPERTIES OF LINEAR MAPPING INDUCED BY LINEAR DESCRIPTOR DIFFERENTIAL EQUATION}\\[5pt]
  Serhiy Zhuk\\
        Cybernetics Faculty\\Taras Shevchenko Kyiv
    National University, Ukraine\\e-mail: serhiy.zhuk@gmail.com
\end{center}

\textbf{Abstract.} In this paper we introduce linear mapping $\oD$ from $\WnF\subset\Ln$ into $\Lm\times\Rm$, induced by linear differential equation $\frac d{dt} F\xt-C(t)\xt=\ft,Fx(t_0)=f_0$. We prove that $\oD$ is closed dense defined mapping for any $m\times n$-matrix $F$. Also adjoint mapping $\oDs$ is constructed and it's domain $\WmF$ is described. 

Some kind of so-called "integration by parts" formula for vectors from $\WnF,\WmF$ is suggested. We obtain  a necessary and sufficient condition for existence of generalized solution of equation $\oD\x=(\f,f_0)$.
 Also we find a sufficient criterion for closureness of the $\R\oD$ in $\Lm\times\Rm$ which is formulated in 
terms of transparent conditions for blocks of matrix $C(t)$. Some examples are supplied to illustrate obtained results.

\subsection*{Introduction}
System of linear differential equations in the form of
\begin{equation}
\label{eq:Ftx}
F(t)\dot x(t)+C(t)x(t)+B(t)f(t)=0
\end{equation}
is called singular or descriptor one. American mathematicians Campbell and Petzold~\cite{cmpbl}  introduced a notion of \textbf{central canonical form} for stationary system~\eqref{eq:Ftx}. Namely, if $\mathrm{det}(\lambda F+C)\ne0$ for any real $\lambda$ then we can transform~\eqref{eq:Ftx} into independent differential and algebraic equations (for sufficiently smooth $\f$ )
$$
 \begin{aligned}
   &\dot x_1(t)= A x_1(t)+K f(t),\,
   x_2(t)=-D f(t)-\sum_{i=1}^{m-1}N^i D v^{(i)}(t)),\,
   \bigl(\begin{smallmatrix}  x_1\\x_2\end{smallmatrix}\bigr)=
   Q^{-1}x(t)
 \end{aligned}
$$
In case of non-constant coefficients in~\eqref{eq:Ftx} russian mathematicians Bojarintsev and Chistjakov suggested a notion of \textbf{left regularization operator}
$$
\Lambda_{*,r}[\frac d{dt}(F(t)x(t))+C(t)x(t)]=\frac d{dt}\xt+\Lambda_{*,r}[C(t)]\xt,
$$ where $\Lambda_{*,r}=\sum_{j=0}^r L_j(t)(\frac d{dt})^i$. In case of constant coefficients existence of central canonical form is equal to left regularization operator existence. In general case conditions of existanse of left regularization operator depends on properties of "prolonged"-system~\cite{chst}. 

Italian mathematician Favini~\cite{fav} studied existence and fundamental solution representation of system~\eqref{eq:Ftx}  in general case where $F(t),C(t)$ supposed to be bounded linear mappings from Banach space $X$ into Banach space $Y$. Their results are based on hypothesis that interval $T:=[a,b]\times[c_1,+\infty)$ consists of regular points of resolvent $(\lambda F(t)+C(t))^{-1}$ which is bounded on $T$.
\looseness=-1
\subsection*{Problem statement}
Papers mentioned above ( and lots of works devoted to singular systems listed in surveys~\cite{muller,surv} ) are based on the hypothesis that some canonical form of~\eqref{eq:Ftx} exists i.e. system~\eqref{eq:Ftx} could be transformed into implicit form. This implies we can use powerfull tool -- theorems about united points for the second type Volterra operators. 

In this paper we use general approach of closed mappings theory combined with regularization methods applied to mapping induced by linear descriptor system 
\begin{equation}
\label{eq:Fx}
\begin{split}
&\frac d{dt} F\xt-C(t)\xt=\ft,\\
&Fx(a)=f_0,
\end{split}
\end{equation} 
Hence we can study some general properties of~\eqref{eq:Fx}  -- existence of solutions, continuous dependence of the solution on the right hand side of operator equation -- without assuming the structure of given system to be canonical. From the other hand operator approach makes it possible to investigate noncasual systems~\cite{muller}.

In~\eqref{eq:Fx} we set $F=\{F_{ij}\}_1^{m,n}$ -- some rectangular matrix, $t\mapsto C(t)$ -- continuous matrix-valued function, $\f$ is some element of squared summable vector-functions space $\Lm:=\LM$, $T<+\infty$, $f_0\in\Rm$. If $F$ is nondegenerate square matrix, then it's easy to see that~\eqref{eq:Fx} has unique totally continuous solution $\x$ and it satisfies Volterra integral equation $$
F\xt=f_0+\int_a^t C(s)x(s)+\f[s]\ds
$$ In general case of rectangular matrix $F$ we define solution of singular initial value problem~\eqref{eq:Fx} as follows. \emph{Set $$F\x[t]=(Fx_1(t),...,Fx_m(t)),Fx_i(t)=\sum_1^n F_{ij}x_j(t)
$$ and let $\WnF$ be a set of all $\x$ from $\Ln$ satisfying $$
F\x \text{\bf is totally continuous and its derivative lies in }\Lm
$$ It's easy to see that $\WnF$ is linear total subset of $\Ln$. For each $\x\in\WnF$ we set $$
\oD\x[t]=(\frac d{dt} F\xt-C(t)\xt,Fx(a))
$$ Now we say that $\x\in\WnF$ is a solution of~\eqref{eq:Fx} if it lies in the solutions domain of operator equation} \begin{equation}
\label{eq:Dxf}
\oD\x=(\f,f_0)
\end{equation}
\emph{The goal of this paper is} investigation of some properties of $\oD$ namely closureness of $\oD$ and conditions for normal solvability of $\oD$. In terms of descriptor systems it can be rewrited as follows: solvability conditions for~\eqref{eq:Fx}, conditions for continuous dependence of solution~\eqref{eq:Fx} on initial condition $f_0$ and perturbation $\f$, approximation of  \eqref{eq:Fx} solution by sequence of functions. 
\subsection*{Closureness of $\oD$ and it's adjoint mapping.}
Now we can introduce 
\begin{thm}\label{t:1}
If\footnote{$\WmF$ is defined in the same way as $\WnF$ with respect to $F'$. } $\x\in\WnF$, $\z\in\WmF$, then
\begin{equation}
\label{eq:ich}
\begin{split}
&\int_a^c(\dF\x[t],\z[t])+(\dFt\z[t],\x[t])\dt=\\
&(F\x[c],F'^+F'\z[c])-(F\x[a],F'^+F'\z[a])
\end{split}
\end{equation}
Moreover, $\oD$ is linear closed dense defined mapping and its adjoint $\oD'$ is given by 
\begin{equation}
\label{eq:Ds}
\begin{split}
&\oD'(\z,z_0)\dfn-\dFt\z-C'\z,\\
&\ov{\oD'}=\{(\z,F'^+F'\z[a]+d):\z\in\WmF,F'\z[c]=0,F'd=0\},
\end{split}
\end{equation} 
\end{thm}
\begin{rem}
We must stress that linear mapping~\footnote{$\Wn$ is a set of totally continuous function in $\Ln$.} $$
\x\mapsto F\dfrac d{dt}\x,\x\in\Wn
$$ is not closed in general case. Really, let's consider case $n=2,t_0=0,T=1$ and set
\begin{equation}
F=\begin{pmatrix}
1&&0\\0&&0
\end{pmatrix}
\end{equation}
If we denote by $t\mapsto k(t)$ Cantor's "dust" function then $v(\cdot):=(0,k(\cdot))\notin\Wn$. Let's set $$
B_n(t)\dfn\sum_0^n k(\frac in){n\choose i}\,t^i\,(1-t)^{n-i}
$$  $v_n(\cdot):=(0,B_n(\cdot))\in\Wn$ and $$
F\dfrac d{dt}v_n(\cdot)\to0,v_n(\cdot)\to v(\cdot)
$$ so $\v\in\Wn$ if $\x\mapsto F\dfrac d{dt}\x$ is closed. \\
On the other hand $\v\in\WnF$ and $\dF\v=(0,0)$.
\end{rem} This remark implies 
\begin{center}
\it In general case $\Wn$ is not Hilbert space with respect to $\dfrac d{dt}F$ graph-norm.
\end{center}
Really, in that case $F\dfrac d{dt}$ would be close on $\Wn$, because $$
\dfrac d{dt}F\x=F\dfrac d{dt}\x,\forall\x\in\Wn
$$
\subsection*{Normal solvability of $\oD$.}
For applications of linear differential equations the range of operator is very important to be closed because it implies continuous dependence of solution on initial conditions and perturbations. Next theorem introduces criterion of \eqref{eq:Fx} pseudosolution existence.
\begin{thm}\label{t:2}
Boundary value problem 
\begin{equation}
\label{eq:FF'e}
\begin{split}
&\dF\x[t]=C(t)\x[t]+\z[t]+\f[t],\\
&\frac d{dt}F'\z[t]=-C'(t)\z[t]+\e^2\x[t],F'\z[c]=0,\\
&F\x[a]-F'^+F'\z[a]+d=f_0,F'd=0
\end{split}
\end{equation}
has unique solution $(\x[\cdot,\e],\z[\cdot,\e],d(\e))$ for any $\e>0$.\\
For given $(\f,f_0)\in\Lm\times\Rm$ descriptor system 
\begin{equation}
  \label{eq:Fx1}
\dF\x[t]=C(t)\x[t]+\ft,F\x[t_0]=f_0
\end{equation} has the
pseudosolution\footnote{We set $\hx\in\WnF$ is the pseudosolution of $\oD\x=(\f,f_0)$ if 
$\|\oD\hx-(\f,f_0)\|^2_2=\min_{\x}\|\oD\hx-(\f,f_0)\|^2_2$. } 
$\hx$ iff $$
\|\x[\cdot,\e]\|_2\le C\text{ while }\e\to0
$$ \end{thm}
\begin{thm}\label{t:3}
Let $$
  F=\bigl(
  \begin{smallmatrix}
    E_r&&0\\0&&0
  \end{smallmatrix}\bigr),
        C(t)\equiv\bigl(\begin{smallmatrix}
    C_1&&C_2\\C_3&&C_4
  \end{smallmatrix}\bigr),f(t)=\bigl(
  \begin{smallmatrix}
    f_1(t)\\f_2(t)
  \end{smallmatrix}\bigr),
        f_0=\bigl(
  \begin{smallmatrix}
    f_1^0\\f_2^0
  \end{smallmatrix}\bigr)
$$
where $E_r$ is identity $r\times r$ matrix, $C_i$ are any matrixes of appropriate dimensions. If\footnote{we set 
$\|F\|_{mod}:=\sum_{i,j}|F_{ij}|$ for any rectangular matrix $F$.} $$
\sup_{1>\e>-1}\|Q(\e)C'_2\|_{mod}<+\infty,Q(\e)\dfn(\e^2E+C'_4C_4)^{-1},
$$ then range of $\oD$ is closed linear manifold.
 \end{thm}
Let's illustrate above theorems by examples.
{\par\noindent\bf Example~1.}
If we set 
\begin{equation*}
   F=\begin{pmatrix}
  1&&0\\0&&0
  \end{pmatrix},
  C(t)\equiv\begin{pmatrix}
  1&&-1\\1&&0
  \end{pmatrix}
\end{equation*}
then $\N\oD=\{0\}$, hence $\cl\R{\oDs}=\Ln$. On the other hand $$
\R{\oDs}=\{
\bigl(\begin{smallmatrix}
  -\dot z_1-z_1-z_2,\\
  z_1
\end{smallmatrix}\bigr),z_1\in\mathbb{W}_2^1([t_0,T]),z_1(T)=0,z_2\in\mathbb{L}_2([t_0,T])\}
$$ according to theorem~\ref{t:1}. So $\R{\oDs}\subset\Ln$ hence range of $\oD$ is not closed. Note that sufficient condition of theorem~\ref{t:3} does not hold in this case because $C'_2Q(\e)=-\e^{-2}$. 

Here~\eqref{eq:FF'e} rewrites as 
\begin{equation}
  \label{eq:xze}
  \begin{split}
    &\dot x_1(t)=x_1(t)+(1+\e^{-2})z_1(t)+f_1(t),\\
    &\dot z_1(t)=-z_1(t)+(1+\e^2)x_1(t)+f_2(t),\\
    &x_1(t_0)-z_1(t_0)=f_{01},z_1(T)=0,x_2(t)=-\e^2 z_1(t,\e)
  \end{split}
\end{equation}
so if we denote by $k(\cdot)$ solution of 
\begin{equation}
  \label{eq:R}
  \dot k(t)=2k(t)+(1+\e^{-2})-(1+\e^2)k^2(t):=U(t,k),k(t_0)=1
\end{equation}
then it's easy to see that 
\begin{equation}
  \label{eq:k+}
  k^-<k(t,\e)<k^+,0<\e<\e_0,t>t_0,
\end{equation}
where $$k^-:=\frac{\e^2-\sqrt{\e^2+3\e^4+\e^6}}{\e^2+\e^4},
k^+:=\frac{\e^2+\sqrt{\e^2+3\e^4+\e^6}}{\e^2+\e^4}$$
So equality $U(t,k)=(1+\e^2)(k-k^-)(k^+-k)$ implies $\dot k(t,\e)>0,t\ge t_0,0<\e<\e_0$, therefore $k(t,\e)\ge k(t_0,\e)>0$ for $t\ge t_0,0<\e<\e_0$.

Let's denote by $q(\cdot,\e)$ solution of $$
q_{tt}(t)-2q_t(t)+(1+\e^{-2})(1+\e^2)q(t)=0,q_t(t_0)=1+\e^2,q(t_0)=1
$$ It's clear that $$
q(t,\e)=e^{\int_{t_0}^t(1+\e^2)k(s,\e)\ds}>0\Rightarrow
q_t(t,\e)\ge0,t\ge t_0,0<\e<\e_0\eqno(*)
$$ If we set $$
\varphi(t,\e)=\frac{e^{t-t_0}}{q(t,\e)}\bigl\{f_1^0+
\int_{t_0}^t \frac{q(\tau,\e)}{e^{\tau-t_0}}f_1(\tau)-
\frac{\dot q(\tau,\e)f_2(\tau)}{e^{\tau-t_0}(1+\e^2)}\mathrm{d}\tau\bigr\}
$$ and $$
z(t,\e)=-\frac{q(t,\e)}{e^t}\int_t^T\frac{e^s}{q(s,\e)}
(f_2(s)+(1+\e^2)\varphi(s,\e))\ds
$$ then it's obvious that $x_1(t,\e)=k(t,\e)z(t,\e)+\varphi(t,\e)$, $x_2(t,\e)=-\e^{-2} z(t,\e)$. Let's set $f_1(t)\equiv0,f_2(t)=-e^{t-t_0},f_1^0=1$. Then $$
\varphi(t,\e)=\frac{\e^2e^{t-t_0}}{(1+\e^2)q(t,\e)}+\frac{e^{t-t_0}}{1+\e^2},
z(t,\e)=-\e^2\frac{q(t,\e)}{e^{t+t_0}}\int_t^T\frac{e^{2s}}{q^2(s,\e)}\ds
$$ and $x(t)=\bigl(\begin{smallmatrix}
  x_1(t)\\x_2(t)
\end{smallmatrix}\bigr),
x_1(t)=-f_2(t)$, $x_2(t)\equiv0$ is unique solution of $\oD\x=(\f,f_0)$.

We'll show that $x_1(\cdot,\e)\to x_1,x_2(\cdot,\e)\to0$ in $\Ln$. $(*)$ implies$$
\frac{e^{t-t_0}}{1+\e^2}<\varphi(t,\e)<\frac{e^{t-t_0}}{1+\e^2}+
\frac{\e^2e^{t-t_0}}{(1+\e^2)q(t_0,\e)},
z(t,\e)\le-\e^2\frac{q(t,\e)}{e^{t+t_0} q^2(t,\e)}\int_t^Te^{2s}\ds,
$$ hence according to~\eqref{eq:k+} and $q(t_0,\e)=1$ we get 
\begin{equation*}
  \begin{split}
    &\int_{t_0}^T(\varphi(t,\e)+f_2(t))^2\dt\le
    \int_{t_0}^T(\frac{e^{t-t_0}}{1+\e^2}+
    \frac{\e^2e^{t-t_0}}{(1+\e^2)}-e^{t-t_0})^2\dt\to0,\e\to0,\\
    &\int_{t_0}^Tk^2(t,\e)z^2(t,\e)\dt\le\int_{t_0}^T(-\e^2k^+
    \frac{e^{2T}-e^{2t}}{2e^{t+t_0}})^2\dt\to0,\e\to0
  \end{split}
\end{equation*} 
therefore $x_1(\cdot,\e)\to-f_2(\cdot)$. One can show 
\begin{equation*}
\begin{split}
&q(t,\e)=\frac{(\e^4+\sqrt{\e^2+3\e^4+\e^6})
e^{\frac{\e^2+\sqrt{\e^2+3\e^4+\e^6}}{\e^2}(t-t_0)}}
{2\sqrt{\e^2+3\e^4+\e^6}}+\\
&\frac{(\sqrt{\e^2+3\e^4+\e^6}-\e^4)
e^{\frac{\e^2-\sqrt{\e^2+3\e^4+\e^6}}{\e^2}(t-t_0)}}
{2\sqrt{\e^2+3\e^4+\e^6}}  
\end{split}
\end{equation*}
hence $\|q(\cdot,\e)\|_2\to+\infty$ while $\e\to0$. On the other hand $$
-\e^{-2} z(t,\e)\le\frac{e^{2T}-e^{2t}}{2e^{t+t_0} q(t,\e)},
$$ therefore $x_2(\cdot,\e)\to0$.
{\par\noindent\bf Example~2.}
If we set \begin{equation*}
   F=\begin{pmatrix}
  -2&&6\\2&&-6
  \end{pmatrix},
  C(t)\equiv\begin{pmatrix}
  1&&-3\\2&&-6
  \end{pmatrix}
\end{equation*}
then~\eqref{eq:Fx} may be rewritten as \begin{equation}
  \label{eq:xze}
  \begin{split}
    &\frac d{dt}(-2x_1+6x_2)(t)=x_1(t)-3x_2(t)+f_1(t),\\
    &\frac d{dt}(2x_1-6x_2)(t)=2x_1(t)-6x_2(t)+f_2(t),\\
    &(-2x_1+6x_2)(t_0)=f^0_1,(2x_1-6x_2)(t_0)=f^0_2
  \end{split}
\end{equation}
One can see that $F_1\dfn LFR=
\bigl(\begin{smallmatrix}
  1&&0\\0&&0
\end{smallmatrix}\bigr),C_0\dfn LC(t)R=
\bigl(\begin{smallmatrix}
  1&&0\\0&&0
\end{smallmatrix}\bigr)
$, where $L=\bigl(\begin{smallmatrix}
  -\frac 13&&\frac 16\\\frac 13&&\frac 13
\end{smallmatrix}\bigr),R=\bigl(\begin{smallmatrix}
  0&&\frac 12\\-\frac 13&&\frac 16
\end{smallmatrix}\bigr)$, so equation~\eqref{eq:xze} is equal to $$
\frac d{dt} F_1R^{-1}\x[t]=C_0R^{-1}\x[t]+Lf(t),F_1R^{-1}\x[t_0]=Lf_0
$$
Note that $\mathrm{det}(\lambda F_1+ C_0)\equiv0$. On the other hand theorem 3 ( $C_2Q(e)\equiv0$ ) implies that range of $y(\cdot)\mapsto\oD_1y(\cdot)=(\frac d{dt} F_1y(\cdot)-C_0y(\cdot),F_1y(t_0))$ is closed. In this case it's simple to verify last sentence. Really, adjoint mapping is defined by rule $$
(\z,z_0)\mapsto
\bigl(\begin{smallmatrix}
  -\dot z_1(t)-z_2(t)\\0
\end{smallmatrix}\bigr),z_2\in\mathbb{L}_2(t_0,T),
z_1\in\mathbb{W}_2^1(t_0,T),z_1(T)=0
$$ so $\R{\oD^*_1}=\mathbb{L}_2(t_0,T)\times\{0\}$ implies $\cl\R{\oD_1}=\R{\oD_1}$. 

\end{document}